\documentclass[10pt,a4paper]{article}
\usepackage{amssymb}
\usepackage{amscd}
\usepackage{latexsym}
\usepackage{amsmath}
\usepackage{amsfonts}
\usepackage{amscd}

\begin{document}
\newtheorem{theorem}{Theorem}[section]
\newtheorem{remark}[theorem]{Remark}
\newtheorem{mtheorem}[theorem]{Main Theorem}
\newtheorem{bbtheo}[theorem]{The Strong Black Box}
\newtheorem{observation}[theorem]{Observation}
\newtheorem{proposition}[theorem]{Proposition}
\newtheorem{lemma}[theorem]{Lemma}
\newtheorem{testlemma}[theorem]{Test Lemma}
\newtheorem{mlemma}[theorem]{Main Lemma}
\newtheorem{note}[theorem]{{\bf Note}}
\newtheorem{steplemma}[theorem]{Step Lemma}
\newtheorem{corollary}[theorem]{Corollary}
\newtheorem{notation}[theorem]{Notation}
\newtheorem{example}[theorem]{Example}
\newtheorem{definition}[theorem]{Definition}

\renewcommand{\labelenumi}{(\roman{enumi})}
\def\Pf{\smallskip\goodbreak{\sl Proof. }}

\def\Fin{\mathop{\rm Fin}\nolimits}
\def\br{\mathop{\rm br}\nolimits}
\def\fin{\mathop{\rm fin}\nolimits}
\def\Ann{\mathop{\rm Ann}\nolimits}
\def\Aut{\mathop{\rm Aut}\nolimits}
\def\End{\mathop{\rm End}\nolimits}
\def\bfb{\mathop{\rm\bf b}\nolimits}
\def\bfi{\mathop{\rm\bf i}\nolimits}
\def\bfj{\mathop{\rm\bf j}\nolimits}
\def\df{{\rm df}}
\def\bfk{\mathop{\rm\bf k}\nolimits}
\def\bEnd{\mathop{\rm\bf End}\nolimits}
\def\iso{\mathop{\rm Iso}\nolimits}
\def\id{\mathop{\rm id}\nolimits}
\def\Ext{\mathop{\rm Ext}\nolimits}
\def\Ines{\mathop{\rm Ines}\nolimits}
\def\Hom{\mathop{\rm Hom}\nolimits}
\def\bHom{\mathop{\rm\bf Hom}\nolimits}
\def\Rk{ R_\k-\mathop{\bf Mod}}
\def\Rn{ R_n-\mathop{\bf Mod}}
\def\map{\mathop{\rm map}\nolimits}
\def\cf{\mathop{\rm cf}\nolimits}
\def\top{\mathop{\rm top}\nolimits}
\def\Ker{\mathop{\rm Ker}\nolimits}
\def\Bext{\mathop{\rm Bext}\nolimits}
\def\Br{\mathop{\rm Br}\nolimits}
\def\dom{\mathop{\rm Dom}\nolimits}
\def\min{\mathop{\rm min}\nolimits}
\def\im{\mathop{\rm Im}\nolimits}
\def\max{\mathop{\rm max}\nolimits}
\def\rk{\mathop{\rm rk}}
\def\Diam{\diamondsuit}
\def\Z{{\mathbb Z}}
\def\Q{{\mathbb Q}}
\def\N{{\mathbb N}}
\def\bQ{{\bf Q}}
\def\bF{{\bf F}}
\def\bX{{\bf X}}
\def\bY{{\bf Y}}
\def\bHom{{\bf Hom}}
\def\bEnd{{\bf End}}
\def\bS{{\mathbb S}}
\def\AA{{\cal A}}
\def\BB{{\cal B}}
\def\CC{{\cal C}}
\def\DD{{\cal D}}
\def\TT{{\cal T}}
\def\FF{{\cal F}}
\def\GG{{\cal G}}
\def\PP{{\cal P}}
\def\SS{{\cal S}}
\def\XX{{\cal X}}
\def\YY{{\cal Y}}
\def\fS{{\mathfrak S}}
\def\fH{{\mathfrak H}}
\def\fU{{\mathfrak U}}
\def\fW{{\mathfrak W}}
\def\fK{{\mathfrak K}}
\def\PT{{\mathfrak{PT}}}
\def\T{{\mathfrak{T}}}
\def\fX{{\mathfrak X}}
\def\fP{{\mathfrak P}}
\def\X{{\mathfrak X}}
\def\Y{{\mathfrak Y}}
\def\F{{\mathfrak F}}
\def\C{{\mathfrak C}}
\def\B{{\mathfrak B}}
\def\J{{\mathfrak J}}
\def\fN{{\mathfrak N}}
\def\fM{{\mathfrak M}}
\def\Fk{{\F_\k}}
\def\bar{\overline }
\def\Bbar{\bar B}
\def\Cbar{\bar C}
\def\Pbar{\bar P}
\def\etabar{\bar \eta}
\def\Tbar{\bar T}
\def\fbar{\bar f}
\def\nubar{\bar \nu}
\def\rhobar{\bar \rho}
\def\Abar{\bar A}
\def\a{\alpha}
\def\b{\beta}
\def\g{\gamma}
\def\w{\omega}
\def\e{\varepsilon}
\def\o{\omega}
\def\va{\varphi}
\def\k{\kappa}
\def\m{\mu}
\def\n{\nu}
\def\r{\rho}
\def\f{\phi}
\def\hv{\widehat\v}
\def\hF{\widehat F}
\def\v{\varphi}
\def\s{\sigma}
\def\l{\lambda}
\def\lo{\lambda^{\aln}}
\def\d{\delta}
\def\z{\zeta}
\def\th{\theta}
\def\t{\tau}
\def\ale{\aleph_1}
\def\aln{\aleph_0}
\def\Cont{2^{\aln}}
\def\nld{{}^{ n \downarrow }\l}
\def\n+1d{{}^{ n+1 \downarrow }\l}
\def\hsupp#1{[[\,#1\,]]}
\def\size#1{\left|\,#1\,\right|}
\def\Binfhat{\widehat {B_{\infty}}}
\def\Zhat{\widehat \Z}
\def\Mhat{\widehat M}
\def\Rhat{\widehat R}
\def\Phat{\widehat P}
\def\Fhat{\widehat F}
\def\fhat{\widehat f}
\def\Ahat{\widehat A}
\def\Chat{\widehat C}
\def\Ghat{\widehat G}
\def\Bhat{\widehat B}
\def\Btilde{\widetilde B}
\def\Ftilde{\widetilde F}
\def\restr{\mathop{\upharpoonright}}
\def\to{\rightarrow}
\def\arr{\longrightarrow}
\def\LA{\langle}
\def\RA{\rangle}
\newcommand{\norm}[1]{\text{$\parallel\! #1 \!\parallel$}}
\newcommand{\supp}[1]{\text{$\left[ \, #1\, \right]$}}
\def\set#1{\left\{\,#1\,\right\}}
\newcommand{\mb}{\mathbf}
\newcommand{\wt}{\widetilde}
\newcommand{\card}[1]{\mbox{$\left| #1 \right|$}}
\newcommand{\union}{\bigcup}
\newcommand{\inters}{\bigcap}
\newcommand{\ER}{{\rm E}}
\def\Proof{{\sl Proof.}\quad}
\def\fine{\ \black\vskip.4truecm}
\def\black{\ {\hbox{\vrule width 4pt height 4pt depth
0pt}}}
\def\fine{\ \black\vskip.4truecm}
\long\def\alert#1{\smallskip\line{\hskip\parindent\vrule%
\vbox{\advance\hsize-2\parindent\hrule\smallskip\parindent.4\parindent%
\narrower\noindent#1\smallskip\hrule}\vrule\hfill}\smallskip}

\title{On Abelian Groups Having Isomorphic \\ Proper Strongly Invariant Subgroups}
\author{Andrey R. Chekhlov \\Faculty of Mathematics and Mechanics, Section of Algebra, \\Tomsk State University, Tomsk 634050, Russia\\{\small e-mails: a.r.che@yandex.ru, cheklov@math.tsu.ru}
\\and\\ Peter V. Danchev \\Institute of Mathematics and Informatics, Section of Algebra, \\Bulgarian Academy of Sciences, Sofia 1113, Bulgaria\\{\small e-mails: danchev@math.bas.bg, pvdanchev@yahoo.com}}
\maketitle

\begin{abstract}{We consider two variants of those Abelian groups with all proper strongly invariant subgroups isomorphic and give an in-depth study of their basic and specific properties in either parallel or contrast to the Abelian groups with all proper fully invariant (respectively, characteristic) subgroups isomorphic, which are studied in details by the current authors in Commun. Algebra (2015) and in J. Commut. Algebra (2023). In addition, we also explore those Abelian groups having at least one proper strongly invariant subgroup isomorphic to the whole group.}
\end{abstract}

\section{Introduction and Definitions}

Throughout the present paper, let all groups into consideration be {\it additively} written and {\it Abelian}. Our notations and terminology from group theory are mainly standard and follow those from \cite{F}, \cite{F1} and \cite{Kap}, respectively. Another useful sources on the explored subject are \cite{CD1,CD2,CD3,CD4} as well. For instance, if $p$ is a prime integer and $G$ is an arbitrary group, $p^nG=\{p^ng~|~ g\in G\}$ denotes the {\it $p^n$-th power subgroup} of $G$ consisting of all elements of $p$-height greater than or equal to $n\in \mathbb{N}$, $G[p^n]=\{g\in G~|~ p^ng=0, n\in \mathbb{N}\}$ denotes the {\it $p^n$-socle} of $G$, and $G_p=\cup_{n<\o} G[p^n]$ denotes the {\it $p$-component} of the {\it torsion part} $tG=\oplus_p G_p$ of $G$.

On the other hand, if $G$ is a torsion-free group and $a\in G$, then let $\chi_G(a)$ denote the {\it characteristic} and let $\t_G(a)$ denote the {\it type} of $a$, respectively. Specifically, the class of equivalence in the set of all characteristics is just called {\it type} and we write $\t$. If $\chi_G(a) \in \t$, then we write $\t_G(a)=\t$, and so $\t(G)=\{\t_G(a) ~ | ~ 0\neq a\in G\}$ is the set of types of all non-zero elements of $G$. The set $G(\t)=\{g\in G ~|~ \t(g)\geqslant \t\}$ forms a pure fully invariant subgroup of the torsion-free group $G$. Recall that a torsion-free group $G$ is called {\it homogeneous} if all its non-zero elements have the same type.

Concerning ring theory, suppose that all rings which we consider are {\it associative} with {\it identity} element. For any ring $R$, the letter $R^{+}$ will denote its {\it additive group}. To simplify the notation and to avoid a risk of confusion, we shall write $\ER (G)$ for the endomorphism ring of $G$ and $\End (G)=\ER (G)^{+}$ for the endomorphism group of $G$.

As usual, a subgroup $F$ of a group $G$ is called {\it fully invariant} (hereafter abbreviated as a {\it fi-subgroup} for simpleness) if $\phi(F)\subseteq F$ for any $\phi\in \ER (G)$, while if $\phi$ is an invertible endomorphism (= an automorphism), then $F$ is called a {\it characteristic} subgroup. Likewise, imitating \cite{GC}, we shall say that a subgroup $S$ of $G$ is {\it strongly invariant}, provided that $\psi(S)\subseteq S$ for any homomorphism $\psi: S\to G$, and in what follows we shall abbreviate it for short by a {\it si-subgroup}. It is well-known that the following relations are fulfilled:

\medskip

\centerline{strongly invariant $\Rightarrow$ fully invariant $\Rightarrow$ characteristic.}

\medskip

Classical examples of important fully invariant subgroups of an arbitrary group $G$ are the defined above subgroups $p^nG$ and $G[p^n]$ for any natural $n$ as well as $tG$ and the maximal divisible subgroup $dG$ of $G$; actually $dG$ is a fully invariant direct summand of $G$ (see, for instance, \cite{F}). To avoid any confusion and misunderstanding, we shall say that a group $G$ has only {\it trivial fully invariant subgroups} if $\{0\}$ and $G$ are the only ones. Same appears for the characteristic and strongly invariant subgroups, respectively.

Let us notice that the si-subgroups were intensively studied in \cite{GC, C1}, respectively, as well as in some other relevant articles cited therewith.

\medskip

Note also that, for all subgroups $F\leq G$, the subgroup $$\mathrm{Hom}(F,G)F=\sum_{\varphi \in \mathrm{Hom}(F,G)}\varphi(F)$$ is the minimal si-subgroup of $G$ containing $F$. In particular, if $F$ is a si-subgroup, then $F=\mathrm{Hom}(F,G)F$. Likewise, we also set $$S_A(B)=\sum_{f\in\mathrm{Hom}(A,B)}\im f$$
to be the $A$-\emph{socle} of $B$.

\medskip

The following key notions, necessary for our successful presentation, were stated in \cite{CD1}.

\medskip

\noindent{\bf Definition 1}. A non-zero group $G$ is said to be an {\it IFI-group} if either it has only trivial fully invariant subgroups, or all its non-trivial fully invariant subgroups are isomorphic otherwise.

\medskip

\noindent{\bf Definition 2}. A non-zero group $G$ is said to be an {\it IC-group} if either it has only trivial characteristic subgroups, or all its non-trivial characteristic subgroups are isomorphic otherwise.

\medskip

Note that Definition 2 implies Definition 1. In other words, any IC-group is an IFI-group; in fact every fully invariant subgroup is characteristic.

\medskip

\noindent{\bf Definition 3}. A non-zero group $G$ is said to be a {\it strongly IFI-group} if either it has only trivial fully invariant subgroups, or all its non-zero fully invariant subgroups are isomorphic otherwise.

\medskip

\noindent{\bf Definition 4}. A non-zero group $G$ is said to be a {\it strongly IC-group} if either it has only trivial characteristic subgroups, or all its non-zero characteristic subgroups are isomorphic otherwise.

\medskip

Notice that Definition 4 implies Definition 3.

On another vein, Definition 4 obviously yields Definition 2, but the converse is false. In fact, in \cite{Kap} was constructed a single non-trivial characteristic subgroup of a $2$-group that are pairwise non-isomorphic, thus giving an example of an IC-group which is surely {\it not} a strongly IC-group; in other words, Definition 4 properly implies Definition 2.

\medskip

We now arrive at our basic tools as follows:

\medskip

\noindent{\bf Definition 5}. A non-zero group $G$ is called an {\it ISI-group} if either it has only trivial strongly invariant subgroups (namely, $\{0\}$ or $G$), or all its non-trivial strongly invariant subgroups are isomorphic otherwise.

\medskip

\noindent{\bf Definition 6}. A non-zero group $G$ is called a {\it strongly ISI-group} if either it has only trivial strongly invariant subgroups, or all its non-zero strongly invariant subgroups are isomorphic otherwise.

\medskip

It is clear that any strongly ISI-group is either a $p$-group for some prime $p$, or is a homogeneous torsion-free group.

\medskip

It is apparent that the next relationships are valid:

\medskip

\centerline{IC-groups $\Rightarrow$ IFI-groups $\Rightarrow$ ISI-groups.}

\medskip

Our objective in the current article is to explore some fundamental and exotic properties of the defined above classes of groups, especially the ISI-groups and the strongly ISI-groups. In addition, we shall investigate even something more, namely the existence of a non-trivial strongly invariant subgroup of a given group, which subgroup is isomorphic to the whole group, calling these groups {\it weakly ISI-groups}. Our motivation to do that is to exhibit and compare certain similarities and discrepancies of these group classes.

The major results established by us are formulated and proved in the next section. Concretely, our work is organized thus: In the next section, we provide our basic statements and their proofs by breaking the results into three parts in conjunction with the definitions stated above and pertained to the three types of the ISI property, namely to the so-termed (strongly, weakly) ISI-groups. In order to prove our results, we shall use some specific ideas and instruments to materialize them (see Theorems~\ref{torsion}, \ref{StrISI} as well as the corresponding lemmas, propositions and examples). We end off our presentation with a series of six questions which, hopefully, will motivate a further intensive study of the explored subject (see Problems 1-6).

\section{Main Theorems and Examples}

For completeness of the exposition and for the reader's convenience, we first and foremost will give a brief retrospection of the most principal results achieved in \cite{CD1} and \cite{CD4}, respectively, concerning IFI-groups and strongly IFI-groups as well as IC-groups and strongly IC-groups, respectively.

\medskip

As usual, the symbol $\oplus_{m} G=G^{({m})}$ will denote the {\it external} direct sum of $m$ copies of the group $G$, where $m$ is some ordinal (finite or infinite).

\medskip

We begin with a retrospection of some of the results obtained in \cite{CD1}.

\begin{theorem}\label{square}
Let $G$ be a $p$-group and let $m\geqslant 2$ be an ordinal. Then $G^{({m})}$ is an IFI-group if, and only if, $G$ is an IC-group.
\end{theorem}

\begin{proposition}\label{tf} Let $G$ be a torsion-free group. Then $G$ is an IFI-group if, and only if, $G$ is a strongly IFI-group.
\end{proposition}

\begin{lemma}\label{need} (a) A fully invariant subgroup of an IFI-group is an IFI-group.

(b) A fully invariant subgroup of a strongly IFI-group is a strongly IFI-group.

\end{lemma}

\begin{proposition} A non-zero IFI-group is either divisible or reduced.
\end{proposition}

\begin{theorem}\label{main} The following two points hold:

(i) A non-zero group $G$ is an IFI-group if, and only if, one of the following holds:

\medskip

$\bullet$ For some prime $p$ either $pG=\{0\}$, or $p^2G=\{0\}$ with $\mathrm{rank}(G)=\mathrm{rank}(pG)$.

\medskip

$\bullet$ $G$ is a homogeneous torsion-free IFI-group of idempotent type.

\medskip

(ii) A non-zero torsion group $G$ is a strongly IFI-group if, and only if, it is an elementary $p$-group for some prime $p$.

\end{theorem}

\begin{proposition}
Every homogeneous fully transitive torsion-free group of idempotent type is an IFI-group.
\end{proposition}

\begin{corollary} A direct summand of a fully transitive torsion-free IFI-group is again a fully transitive IFI-group.
\end{corollary}

We now continue with some of the most important corresponding assertions from \cite{CD4}.

\begin{proposition} Every torsion IFI-group is an IC-group.
\end{proposition}

\begin{theorem} There exists an IFI-group which is not an IC-group.
\end{theorem}

\subsection{ISI-groups}

We start our work here with some elementary but useful observations:

\begin{proposition}\label{only} Let $G$ be a non si-simple group. Then $G$ is an ISI-group if, and only if, $G$ has a single non-trivial si-subgroup $H$.
\end{proposition}

\Pf "$\Rightarrow$". Assume that both $H$ and $F$ are non-trivial si-subgroups. Thus, $H\cong F$, so $F\leq\mathrm{Hom}(H,G)H$ and hence $F\leq H$. Analogously, we derive that $H\leq F$, i.e., $H=F$, as required.

"$\Leftarrow$". It is evident, so we omit the details.
\fine

The next Theorem is quite similar to \cite[Theorem 2.5]{CD1}, but however is somewhat its reminiscent.

\begin{theorem}\label{torsion}
A non-zero torsion group is an ISI-group if, and only if, for some prime $p$, either $pG=\{0\}$ or $p^2G=\{0\}$.
\end{theorem}

\Pf Suppose first that $G$ is torsion, that is, $G=tG$. If $G=G[p]$, the assertion follows automatically. So, assume next that $G\neq G[p]$. We now assert that $G=G[p^2]$. If $G\neq G[p^2]$, then $G[p]\cong G[p^2]$ which is an absurd, so that the claim is sustained.

Reciprocally, if $G$ is an elementary $p$-group, then it contains only trivial si-subgroups, and thus we are done. So, let $G$ be $p^2$-pounded. It is well known in this case that the only proper si-subgroups of $G$ is $G[p]$, as needed.
\fine

In some cases the subgroup $H$ from Proposition~\ref{only} is si-simple (which is the case not only in Theorem~\ref{torsion}).

\begin{proposition}\label{fids}
If $G=A\oplus B$, where $A\neq\{0\}$ is fully invariant in $G$, then $G$ is an ISI-group if, and only if, $A$ and $B$ are such si-simple groups that $\mathrm{Hom}(B,A)\neq\{0\}$.
\end{proposition}

\Pf "$\Rightarrow$". If we assume in a way of contradiction that $A$ is not a si-simple group, then it will contain an si-subgroup, say $\{0\}\neq F\neq A$. But then both $A$ and $F$ are simultaneously si-subgroups in $G$, so that, by assumption, we can deduce that $F\cong A$. In fact, if $S=\mathrm{Hom}(F,B)F$, then $A\oplus S$ is a si-subgroup in $G$, whence $A\cong A\oplus S$, but thus contradicting $\mathrm{Hom}(A,B)=\{0\}$. So, with no harm of generality, we may assume that $A\leq F$, i.e., $F=A$, which is the desired contradiction.

Also, if $\mathrm{Hom}(B,A)=\{0\}$, then both $A$ and $B$ are non-isomorphic si-subgroups. If, however, $\{0\}\neq C\neq B$ is a si-subgroup, then $A\oplus C\cong A$ and, therefore, $\mathrm{Hom}(A,B)\neq \{0\}$, again
a contradiction, whence $B$ is si-simple, as wanted.

"$\Leftarrow$". Let $\{0\}\neq H\leq G$ be a si-subgroup and let $\pi: G\to B$ be the corresponding projection. If $H\leq A$, then $H$ is just a si-subgroup in $A$, and so $H\cong A$. Assuming now that $H\nleq A$, then
the relations $\{0\}\nleq \pi(H)\leq B$ hold. So, $$B\leq\mathrm{Hom}(\pi(H),B)\pi(H)\leq \mathrm{Hom}(\pi(H),G)\pi(H)\leq H$$ and, consequently, $H=B\oplus(H\cap A)$, where $H\cap A$ is a si-subgroup in $A$. Finally, $H\cap A=A$
and thus $H=G$, as pursued.
\fine

It follows from Proposition~\ref{fids} that the sufficiency in the following two statements are directly true.

\begin{theorem}\label{StrISI} The mixed group $G$ is an ISI-group if, and only if, $G=T\oplus R$, where $T$ is an elementary $p$-group for some $p$ and $R$ is such a si-simple torsion-free group that $pR\neq R$.
\end{theorem}

\Pf "$\Rightarrow$". It is clear that the torsion part $tG$ of such a group $G$ is an elementary $p$-group, so one can decompose $G=T\oplus R$, where $T=tG$ and $R$ is torsion-free. If $F\neq \{0\}$ is a si-subgroup in $R$ and $F\neq R$, then $T$ and $T\oplus F$ are apparently non-isomorphic si-subgroups of $G$, so that $R$ is a si-simple group, as claimed.
\fine

\begin{proposition}\label{non-reduced} The non-reduced group $G$ is an ISI-group if, and only if, $G=D\oplus R$, where $D$ is a torsion-free divisible group and $R$ is a si-simple torsion-free group.
\end{proposition}

\Pf "$\Rightarrow$". Clearly, the divisible part of $G$ is torsion-free (by following the same arguments as in Theorem~\ref{StrISI}) and, moreover, we may apply the same Theorem~\ref{StrISI} to get that the group $R$ is a si-simple, as asserted.
\fine

We shall say that a torsion-free group $G$ is {\it strongly irreducible} (or, shortly, {\it s-irreducible}) if $G$ does not have proper pure si-subgroups.

It is very clear that if $G$ is a torsion-free group and $nG=G$ for some $n\in\mathbb{N}$, then $nH=H$ for each fi-subgroup $H$ of $G$. This follows immediately from the fact that $n^{-1}\cdot 1_G\in\mathrm{E}(G)$.

\medskip

The next two examples are crucial for our considerations.

\begin{example}\label{simplehom} There exists a non si-simple homogeneous torsion-free group.
\end{example}

\Pf Let $A$ be a module on the ring $\widehat{\mathbb{Z}}_p$ consisting of $p$-adic integers such that, as a group, $A$ is reduced torsion-free, and let $B$ be a separable homogeneous torsion-free group of type $\t=\t(\widehat{\mathbb{Z}}^+_p)$. Hence, $G=A\oplus B$ is homogeneous and also it is readily checked that $G$ is an ISI-group which is not si-simple, because $A$ is a si-direct summand in $G$.
\fine

\begin{example} There exists an s-irreducible group that is not an irreducible group.
\end{example}

\Pf Let $G$ be the group from \cite[\S~88, Exersize~6]{F1}. Then, $\mathrm{E}(G)=\mathbb{Z}$ and the rank of $G$ can be chosen to be a free natural number $n\geq 2$. So, every subgroup of $G$ is fully invariant and, consequently, $G$ is not irreducible. Likewise, each subgroup $A$ in $G$ of rank $1\leq rank(A)\leq n-1$ is a free group, so that $\mathrm{Hom}(A,G)A=G$. But, if the subgroup $B$ is of rank $n$, then one inspects that $B$ is essential in $G$, whence $F=\mathrm{Hom}(A,G)A=G$ is also essential in $G$ and thus $F_*=G$. Therefore, $G$ is s-irreducible. Note also that since all pure subgroups in $G$ of rank $\leq n-1$ are free, they are not si-subgroups, thus answering in the negative question (1), posed in \cite{GC}.
\fine

The following technicalities are worthwhile for further applications.

\begin{proposition}\label{pG} If $G$ is a torsion-free ISI-group and $\{0\}\neq H\lneqq G$ is a si-subgroup, then $H$ is $q$-pure in $G$ for all prime numbers $q$ except, maybe, one prime number $p$. Moreover, if $H$ is not $p$-pure in $G$ for some prime $p$, then $pG\lneqq H$.
\end{proposition}

\Pf Assume that $pg\in H$ for some prime $p$, where $g\in G\setminus H$. Letting
$$F=\mathrm{Hom}(\langle H,g\rangle,G)\langle H,g\rangle,$$
then $F$ is a si-subgroup in $G$, so either $F\cong H$ or $F=G$. Besides, the condition $F\cong H$ yields $H=F$ that contradicts the fact that
$g\in G\setminus H$. So, it follows that $F=G$.

Furthermore, each element $y\in G$ can be written in the form
$$y=f_1(x_1+l_1g)+\dots+f_n(x_n+l_ng)$$
for some $f_i\in\mathrm{Hom}(\langle H,g\rangle,G)$, $x_i\in H$ and $l_i\in\mathrm{Z}$; $i=1,\dots,n$. That is why,
$$py=f_1(px_1+pl_1g)+\dots+f_n(px_n+pl_ng)\in H$$
since $px_i+pl_ig\in H$. So, $pG\leq H$, where $pG\neq H$ since $pG\cong G$ and $H\neq G$. If, however, $H$ is not $q$-pure for some prime $q\neq p$, then similarly $qG\leq H$, which implies that $H=G$, so $H$ is indeed $q$-pure for all primes $q\neq p$.
\fine

\begin{remark} Let us note that Proposition~\ref{non-reduced} and Example~\ref{simplehom} manifestly give valuable examples of proper pure si-subgroups in torsion-free ISI-groups. Also, note that the si-subgroup $H$ in Proposition~\ref{pG} has the following property: $$\mathrm{Hom}(A,G)A=H,$$ provided $$\{0\}\neq\mathrm{Hom}(A,G)A\neq G$$ for every group $A$.
\end{remark}

The next assertion is rather curious by comparing with the inheritance of the direct summand property by the IFI and IS groups, respectively.

\begin{lemma}\label{summand} A direct summand of an ISI-group $G$ is again an ISI-group.
\end{lemma}

\Pf Write $G=A\oplus B$ and given $S\leq A$ is a si-subgroup in $A$. If we set $F=\mathrm{Hom}(S,B)S$, then $S\oplus F$ is obviously a si-subgroup in $G$. However, if $S\oplus F\neq G$, then according to Proposition~\ref{only} we conclude that $S\oplus F$ is the only non-trivial si-subgroup.
Hence, it follows that $A$ is an ISI-group, as formulated.
\fine

We also notice that it follows immediately from Theorem~\ref{StrISI} and Lemma~\ref{summand} that the following is true:

\begin{proposition} The non-zero divisible group $D$ is an ISI-group if, and only if, $D$ is torsion-free.
\end{proposition}

We are now ready to give our two desired examples.

\begin{example}\label{pure} Pure si-subgroups of an ISI-group $G$ are not necessarily direct summands of $G$.
\end{example}

\Pf Let $A$ and $B$ be both torsion-free groups of rank $1$ with $\t(A)<\t(B)=\t$, and let
$$G=\langle A\oplus B,p^{-1}(a+b)\,|\,a\in A\setminus pA,b\in B\setminus pB\rangle.$$
Then, $B$ is pure in $G$ and $B=G(\t)$, so $B$ is a si-subgroup in $G$. Suppose that $\pi: A\oplus B\to A$ is a projection. If $H\neq B$ is a si-subgroup of $G$, then one deduces that $H\nleqslant B$ since $B$ is si-simple, and hence $\pi(pH)\neq\{0\}$. Thus, it follows that $A\leq H$ and so $H=G$, as expected.
\fine

Notice that an example of si-subgroups in non-homogeneous torsion-free groups is given in the previous Example~\ref{pure}, and an example of pure si-subgroups in homogeneous torsion-free groups is given in Example~\ref{simplehom}.

\begin{example} There is a homogeneous torsion-free group $A$ of rank $r\geq 2$ with a si-subgroup $H$
as exhibited in Proposition~\ref{pG}.
\end{example}

\Pf Let $A$ be the group as constructed in \cite[Example~5 from \S88]{F1}, that is,
$$A=\langle a_1,\dots,a_r,x_1,x_2,\dots\rangle,$$
where $x_n=p^{-n}(a_1+\pi_{2n}a_2+\dots+\pi_{rn}a_r)$, $\pi_{in}=s_{i0}+s_{i1}p+\dots+s_{i,n-1}p^{n-1}$
is the $n-1$ partial sum of the $p$-adic reversible number $\pi_i$ ($i=2,\dots,r$), $\pi_i$ are algebraically independent on $\mathbb{Q}$ and $\pi_1=1$. It is principally known that $\mathrm{E}(A)=\mathbb{Z}$, the subgroup
$\langle a_2,\dots,a_r\rangle$ is pure in $A$ and all of the elements from $A\setminus A_0$ are of the form
$kx_n+k_2a_2+\dots+k_ra_r$, where $k,k_2,\dots,k_r\in\mathbb{Z}$, $p\nmid k$ and $n$ is some natural number.
Letting
$$H=\langle a_1,a_2\dots,a_r,px_1,px_2,\dots\rangle,$$
it is clear then that $pA\leq H$ and $H\neq A$.

We now manage to show that $\mathrm{Hom}(H,A)=\mathrm{Z}$. To that goal, we shall use some techniques from \cite[Example 5 of \S88]{F1}. And so, choose $0\neq\eta\in\mathrm{Hom}(H,A)$ and define the homomorphism $\eta$ thus:
$$\eta a_i=\sum_{j=1}^rt_{ij}a_j, ~ t_{ij}\in\mathbb{Z}.$$
Notice that, if required, multiply $\eta$ on some $m\in\mathbb{N}$ so that the new map will send $A_0=\langle a_1,\dots,a_r\rangle$ into itself.
Therefore,
$$\eta (px_n)=p^{1-n}\sum_i\pi_{in}\eta a_i=
p^{1-n}\sum_j\Bigl(\sum_i\pi_{in}t_{ij}\Bigr)a_j=
k_npx_n+l_{2n}a_2+\dots+l_{rn}a_r$$
for some
$k_n,l_{2n},\dots,l_{rn}\in\mathbb{Z}$.
We also have
$$\sum_{i=1}^r\pi_{in}t_{i1}=k_n
\  \mathrm{and}\
p^{1-n}\sum_{j=2}^r\Bigr[\sum_{i=1}^r\pi_{in}t_{ij}-k_n\pi_{jn}\Bigr]a_j\in
\langle a_2,\dots,a_r\rangle.$$
Furthermore, one readily can verify that the coefficients in the square brackets are divided by $p^{n-1}$.
Also, under taking $n\rightarrow\infty$, for each index $j=2,\dots,r$ we obtain that $\sum_{i=1}^r\pi_it_{ij}-\kappa\pi_j=0$, where $\sum_{i=1}^r\pi_it_{i1}=\kappa$.

Now, in view of the algebraic independence of the $\pi_i$'s, it follows that $t_{jj}=t_{11}$ and $t_{ij}=0$
assuming $i\neq j$. Thereby, the homomorphism $\eta$ acts as multiplying on the integer $t_{11}$. Since
$\chi(a_i)=(0,0,\dots)$, there is no problem to take $m>1$. So, we infer that $\mathrm{Hom}(H,A)=\mathbb{Z}$, as asked for, and hence $H$ is a si-subgroup in $A$, as wanted.
\fine

The following technicality is pretty simple but useful.

\begin{lemma} (1) If $G$ is a torsion-free group and $H\leq G$ is a si-subgroup, then $H_*$ also is a si-subgroup
in $G$.

(2) A direct summand of s-irreducible (si-simple) group $G$ is again an s-irreducible (si-simple).
\end{lemma}

\Pf (1) Letting $f\in\mathrm{Hom}(H_*,G)$ and $x\in H_*$, we have $nx\in H$ for some $n\in\mathbb{N}$, so that $f(nx)=nf(x)\in H$ whence $f(x)\in H_*$.

(2) Write $G=A\oplus B$ and assume that $H\leq A$ is a si-subgroup in $A$ and put $F=\mathrm{Hom}(H,B)H$. Thus,
$H\oplus F$ is a si-subgroup in $G$. If, however, $H\neq A$, then one has that $H\oplus F\neq G$ that contradicts the condition that the group $G$ is si-simple.

Similarly, if $H_*\neq A$, then we will have that $H_*\oplus F_*\neq G$ which contradicts the condition that the group $G$ is s-irreducible.
\fine

The next two statements are worthy of noticing.

\begin{proposition}\label{s-irr} Let $G$ be a torsion-free s-irreducible group which is not a si-simple ISI-group. Then its non-trivial si-subgroup $H$ is also s-irreducible.
\end{proposition}

\Pf Owing to Proposition~\ref{pG}, we derive that $pG< H$ for some prime $p$. Let $S\neq\{0\}$ be a pure si-subgroup in $H$ and $F=\mathrm{Hom}(S,G)S$. Then, one sees that either $F=H$ or $F=G$. If, for a moment, $F=H$, then
$$H=\mathrm{Hom}(S,G)S=\mathrm{Hom}(S,H)S=S.$$
Assume now that $F=G$. Therefore,
$$pG=\sum_{f\in\mathrm{Hom}(S,G)}pf(S)=\sum_{f\in\mathrm{Hom}(S,H)}pf(S)$$
since all $pf(S)\leq H$. Further, $pG\leq S$ as $\sum_{f\in\mathrm{Hom}(S,H)}f(S)=S$. In particular, $pH\leq S$ and so $S=H$ in view of purity of
$S$ in $H$, as promised.
\fine

The next assertion is closely relevant to Example~\ref{pure}.

\begin{proposition}
(1) The pure si-subgroup of a direct sum $G$ of s-irreducible torsion-free groups coincides with some fully invariant direct summand of $G$.

(2) The direct sum of s-irreducible torsion-free groups $A_i$, where $i\in I$, $|I|\geq 2$, is an s-irreducible group if, and only if, $\mathrm{Hom}(A_i,A_j)\neq\{0\}$ for all indices $i,j\in I$.

(3) The direct sum $G$ of si-simple groups $A_i$, where $i\in I$, $|I|\geq 2$, is a si-simple group if, and only if,
$\mathrm{Hom}(A_i,A_j)\neq\{0\}$ for all indexes $i,j\in I$.
\end{proposition}

\Pf (1) Let us write $G=\bigoplus_{i\in I}A_i$, where all $A_i$ are s-irreducible groups and $H\leq G$ is a si-subgroup. Then, $H=\bigoplus_{i\in I}(H\cap A_i)$. However, since $H_*=\bigoplus_{i\in I}(H\cap A_i)_*$, we have that $H_*=\bigoplus_{j\in J} A_j$, where
$J\subseteq I=\{j\in I\,|\,H\cap A_j\neq\{0\}\}$, as required.

(2) "$\Rightarrow$". Assume that $\mathrm{Hom}(A_i,A_j)=\{0\}$ for some $i\neq j$. Setting $$S=\{s\in I\,|\,\mathrm{Hom}(A_i,A_s)\neq\{0\}\},$$
we observe that $S\neq\emptyset$, because $i\in S$ and $S\neq I$ since $j\in I\setminus S$, so $\bigoplus_{s\in S}A_s$ is a proper pure si-subgroup (i.e., $\neq \{0\},G$) that contradicts the s-irreducibility of $G$.

"$\Leftarrow$". Given $H$ is a pure si-subgroup in $G$. We thus may write $H=\bigoplus_{i\in I}(H\cap A_i)$ and
$H\cap A_i=A_i$ with $H\cap A_i\neq\{0\}$ which holds in view of s-irreducibility of $A_i$ for each index $i$. So, one follows that $H=G$ as $\mathrm{Hom}(A_i,A_j)\neq\{0\}$ for all $i,j\in I$.

(3) The proof is similar to that in (2).
\fine

Furthermore, in virtue of Proposition~\ref{fids}, in the next statement the fulfillment of the inequalities $\mathrm{Hom}(A,B)\neq\{0\}$ and $\mathrm{Hom}(B,A)\neq\{0\}$ is guaranteed.

\begin{proposition} Let $A$, $B$ be ISI-groups, where $\mathrm{Hom}(A,B)\neq\{0\}$, $\mathrm{Hom}(B,A)\neq\{0\}$
and $\{0\}\neq H\leq A$, $\{0\}\neq F\leq B$ are si-subgroups of $A$, $B$, respectively (in particular, $H=A$, $F=B$
provided $A$, $B$ are si-simple groups). Then the group $G=A\oplus B$ is an ISI-group if, and only if, $\mathrm{Hom}(H,B)H=F$ and $\mathrm{Hom}(F,A)F=H$, as moreover $\mathrm{Hom}(A,B)A=B$ and $\mathrm{Hom}(B,A)B=A$
provided both $A$, $B$ are not si-simple.
\end{proposition}

\Pf "$\Rightarrow$". If, for a moment, $$U=\mathrm{Hom}(H,B)H\neq F,$$ then $H\oplus U$ is a si-subgroup in $G$ with
$H\oplus U\neq H\oplus F$ that contradicts the uniqueness of the si-subgroup of an ISI-group. Thus, hereafter, if both $A$ and $B$ are not si-simple, i.e., $H\neq A$ and $F\neq B$ and, for instance, $\mathrm{Hom}(A,B)A=F$, then we have for the si-subgroup $A\oplus F\neq H\oplus F$.

"$\Leftarrow$". Assume that $\{0\}\neq U, V\leq G$ are some si-subgroups in $G$. It suffices to show that $U=V$. To that goal, we have $$U=(U\cap A)\oplus(U\cap B),V=(V\cap A)\oplus(V\cap B).$$ On this vein, since $U\cap A$, $V\cap A$ are si-subgroups in $A$ as well as $U\cap B$, $V\cap B$ are si-subgroups in $B$, and $A$,$B$ are ISI-groups, one has that $U\cap A=V\cap A$ and $U\cap B=V\cap B$ hold, provided that these subgroups are not trivial.

Note that the condition $U\cap A=\{0\}$ is impossible. In fact, if we assume in a way of contradiction that $U\cap A=\{0\}$, then $U\cap B\neq\{0\}$ since $U\neq\{0\}$. If, however, $U\cap B=F$, then
$$\{0\}\neq H=\mathrm{Hom}(F,A)F\leq U\cap A,$$ thus contradicting the equality $U\cap A=\{0\}$. But, if $U\cap B=B$, then $U\cap A=\{0\}$, thus contradicting the condition $\mathrm{Hom}(B,A)\neq\{0\}$. Analogously, one obtains that $U\cap B\neq\{0\}$, $V\cap A\neq\{0\}$, $V\cap B\neq\{0\}$.

Assume now that $U\cap A=A$. Then, the condition $\mathrm{Hom}(A,B)\neq\{0\}$ forces that $U\cap B\neq\{0\}$. So, we have $U=A\oplus F$, where $F=U\cap B\neq B$, because $U\neq G$. In particular, $B$ is not si-simple.

Next, we consider the first case when $A$ is si-simple. To this purpose, assume that $V\cap B=B$. Since $\mathrm{Hom}(B,A)\neq\{0\}$, it follows that $V\cap A\neq\{0\}$ and so $V\cap A=A$ taking into account that $A$ is si-simple. Therefore, it must be that $V=G$, thus contradicting the given assumption on $V$. So, $V\cap B=F\neq B$
and since $V\cap A\neq\{0\}$, we obtain that $V\cap A=A$, i.e., $V=A\oplus F$. Consequently, $U=V$ holds in the case
when $A$ is si-simple.

Finally, assume $A$ is not si-simple. If, however, $B$ is si-simple, then we deduce that $U\cap B=B$ and so $U=G$ follows, a contradiction. So, the condition $U\cap A=A$ follows, provided that $B$ is not si-simple. By assumptions, we know that $\mathrm{Hom}(A,B)A=B$ and $\mathrm{Hom}(B,A)B=A$. Thus, the equality $U\cap A=A$ assures that $U\cap B=B$, i.e., $U=G$, thus contradicting $U\neq G$, as required.
\fine

As an immediate consequence, we extract the following:

\begin{corollary}\label{cartesian} If $G$ is either a si-simple, s-irreducible torsion-free group or an ISI-group, then, for every ordinal $m>1$, the group $G^{(m)}$ retains the same properties.
\end{corollary}

Note that some other interesting results for si-subgroups and irreducible groups are proved in \cite{C2,C1}, respectively. We now will give some additions to that as follows. Let us recall that, for a group $G$, the notation $\tau(G)$ means the set consisting of all types of the non-zero elements of $G$.

\begin{proposition}\label{01new} If $G$ is a non-homogeneous torsion-free ISI-group, then $\tau(G)=\{t_1,t_2\}$, where $t_1<t_2$ and the non-trivial si-subgroup of $G$ is exactly $G(t_2)$.
\end{proposition}

\Pf We clearly have $\{0\}\neq G(t_2)\neq G$ for some type $t_2$. Since the si-subgroup is unique (see Proposition~\ref{only}), the group $G(t_2)$
is then homogeneous, i.e., the type $t_2$ is maximal. Letting $t_1\in \tau(G)\setminus\{t_2\}$, it follows from $G(t_1)\ncong G(t_2)$ that $G(t_1)=G$. Consequently, we deduce that $t_1<t_2$ and $\tau(G)=\{t_1,t_2\}$, as stated.
\fine

Note that the si-subgroup $H=G(t_2)$ from Proposition~\ref{01new} is manifestly s-irreducible. Indeed, if $\{0\}\neq F<H$ is a pure si-subgroup in $H$, then we can plainly inspect that $\{0\}\neq\mathrm{Hom}(F,G)F\leq H$, so that $F=H$, as claimed.

\medskip

We are now recording the following consequence.

\begin{corollary} If the torsion-free ISI-group $G$ has a non-$p$-pure si-subgroup $H$, then $G$ is homogeneous and, if $H$ is not si-simple, then $pG\leq F$ for every non-trivial si-subgroup $F$ of $H$.
\end{corollary}

\Pf Arguing as above, we are aware that
$$\sum_{\alpha\in\mathrm{Hom}(F,G)}\im\,\alpha=G,$$
whence
$$pG=\sum_{\alpha\in\mathrm{Hom}(F,G)}\im\,(p\alpha)\leq F$$
as
$p\alpha\in\mathrm{Hom}(F,H)$ by exploiting Proposition~\ref{pG}.
\fine

Now, we are able to show the validity of the following two constructions.

\begin{example} The next two points are true:

(i) For every type $\t\neq (\infty,\infty,\dots)$, there exists a homogeneous torsion-free group $G$ of type $\t$ such that, for each natural number $r>1$, $G$ has a si-subgroup of rank $\frac{(r+2)(r-1)}{2}$.

(ii) There exists a torsion-free group $G$ such that, for each natural number $r\geq 2$, $G$ possesses a si-subgroup of rank $2r$.
\end{example}

\Pf (i) Let $G_r$ be the group from \cite[\S 88, Example~5]{F1} having rank $r$ and set $G=\bigoplus_{r>1} G_r$. Then, the so-constructed group $G$ will have type $(0,0,\dots)$. By virtue of \cite[\S 88, Exersize~6]{F1}, any subgroup in $G_r$ of rank $\leq r-1$ is a free group, and hence it follows that the following equality is true
$$\mathrm{Hom}(\bigoplus_{2\leq i\leq r}G_r,\bigoplus_{j>r}G_j)=\{0\},$$
because $G_r$ does not have unfree direct summands, i.e., each direct component $\bigoplus_{2\leq i\leq r}G_r$ is a si-subgroup of $G$.

If now $R$ is a torsion-free group of rank $1$ and type $\t\neq (\infty,\infty,\dots)$, then one readily verifies that the standard tensor product $G\otimes R$ over $\mathbb{Z}$ is a homogeneous group of type $\t$, equipped with the same properties as these of $G$, as needed.

(ii) Let $G_r$ be the group of rank $2r$ from \cite[\S 88, Exersize~7]{F1} and put $G=\bigoplus_{r\geq 2} G_r$. Then, since every subgroup in $G_r$ of rank $\leq r$ is free and all factor-groups of $G_r$ of rank $\leq r$ are divisible, one easily checks that each $G_r$ is a si-subgroup in $G$, as required.
\fine

We shall now discover when ISI-groups are IFI-groups by arranging to prove the following (as mentioned in the introductory section, the converse is always true, because si-subgroups are fi-subgroups).

\begin{proposition}\label{conver} (i) A torsion ISI-group $G$ is an IFI-group if, and only if, $G$ is a $p$-group for some prime $p$ having the property $p^2G=\{0\}$ and, moreover, if $pG\neq\{0\}$, then $\mathrm{rank}(G)=\mathrm{rank}(pG)$.

(ii) A mixed ISI-group need not be an IFI-group.

(iii) A torsion-free IFI-group is a si-simple group of idempotent type.
\end{proposition}

\Pf (i) It follows automatically from Theorems~\ref{main} and \ref{torsion}.

(ii) It follows directly from a combination of Theorem~\ref{main} (since there not exists a mixed IFI-group) and Theorem~\ref{StrISI} in which mixed ISI-groups are described.

(iii) It follows immediately from the fact in Proposition~\ref{tf} that torsion-free IFI-groups are strongly IFI-groups, as in this case every si-subgroup is isomorphic to the whole group and, consequently, all torsion-free IFI-groups are themselves si-simple.
\fine

Note that si-simple groups are always ISI-groups since they have only trivial si-subgroups. Also, there are si-simple groups that are {\it not} IFI-groups; indeed, all separable torsion-free groups of a non-idempotent type are such groups. However, a question which immediately arises is whether all si-simple torsion-free groups of idempotent type are IFI-groups?

\subsection{Strongly ISI-groups}\label{strong}

We continue here with more details than the comments listed after Definition 6.

\begin{lemma}\label{01} The following conditions are equivalent:

(a) The group $G$ is a strongly ISI-group;

(b) The equality $\mathrm{Hom}(F,G)F=G$ holds for every subgroup $\{0\}\neq F\leq G$;

(c) The group $G$ is a si-simple group.
\end{lemma}

\Pf (a) $\Rightarrow$ (c). Indeed, if $H\neq \{0\}$ is a si-subgroup, then $H\cong G$, and if $f: H\to G$ is an isomorphism, then $f(H)=G\leq H$, i.e., $H=G$, as needed.

The relationships (b) $\Leftrightarrow$ (c) and (c) $\Rightarrow$ (a) are obvious, so we omit the details in proving them.
\fine

It is clear that a si-simple torsion-free group is a homogeneous group and that every fi-simple group is si-simple. Moreover, in view of \cite[Proposition~24]{GC}, each fi-simple group is either $p$-elementary for some prime $p$ or is divisible torsion-free.

\medskip

Some other examples of si-simple groups are the following ones:

\medskip

(1) Fully transitive homogeneous torsion-free groups of idempotent type are si-simple (see, e.g., \cite[Proposition 26]{GC}).

(2) For any type $\t$, there exists a torsion-free si-simple group of homogeneous type $\t$. In fact, as such a group, it is possible to take any separable homogeneous torsion-free group $G$ (see, for instance, \cite[Proposition 1(1)]{C1}).

(3) A non-zero torsion group $G$ is a si-simple group if, and only if, it is a non-zero elementary $p$-group for some prime $p$. This claim follows at once from \cite[Proposition 24]{GC}.

\medskip

Furthermore, we give here some elementary but helpful observations:

\medskip

$\bullet$ A torsion-free IFI-group is always a strongly ISI-group.

\medskip

Indeed, we just need to combine Proposition~\ref{conver} (iii) and Lemma~\ref{01} (c).

\medskip

$\bullet$ A non-zero torsion group is a strongly ISI-group if, and only if, it is a non-zero elementary $p$-group for some prime $p$.

\medskip

Indeed, it suffices just to combine point (3) quoted above with Lemma~\ref{01} (c).

\subsection{Weakly ISI-groups}

Before starting our work, it is worth to notice that in \cite{GN1}, \cite{GN2} and \cite{GN3} were examined those Abelian groups, both torsion and torsion-free, having isomorphic proper fully invariant subgroup. Likewise, in \cite{Fo} were investigated in detail those Abelian groups in which all subgroups of infinite index are free. 

In this subsection, we will initiate the examination of the case when we have a proper strongly invariant subgroup isomorphic to the whole group. For simpleness of the exposition, we shall call these groups just {\it weakly ISI-groups}. However, unfortunately, this class of groups is definitely {\it not} so interesting, because each si-subgroup which is a weakly ISI-group coincides with the whole group.

\section{Concluding Discussion and Open Problems}

We conclude this final section with some valuable comments on the obtained results. In fact, summarizing all of what we have established so far, we can just say that the properties of ISI-groups are totally different from these of ISI-groups and IC-groups which were completely described in \cite{CD1} and \cite{CD4}, respectively. The primary reason for this discrepancy is that there is no abundance of so many si-subgroups (compare also with \cite{GC}).

\medskip

We close our work with six questions of interest and importance.

\medskip

\noindent{\bf Problem 1.} Does there exist s-irreducible that are {\it not} ISI-groups?

\medskip

\noindent{\bf Problem 2.} Is it true that si-subgroups of ISI-groups are also ISI-groups?

\medskip

\noindent{\bf Problem 3.} Does it follow that an s-irreducible torsion-free ISI-group is si-simple?

\medskip

Knowing with the aid of \cite[Example 2.6]{CD1} that there exists a homogeneous torsion-free group of idempotent type and rank strictly greater that $1$ which is {\it not} an IFI-group, it is reasonably logical to state the following.

\medskip

\noindent{\bf Problem 4.} Decide when a homogeneous torsion-free group of idempotent type and rank $2$ is an ISI-group?

\medskip

We know that each strongly invariant subgroup is always fully invariant, and thus every IFI-group is an ISI-group, but the reverse implication fails. So, in regard to Proposition~\ref{conver} and the first bullet from Subsection~\ref{strong}, we may ask the following.

\medskip

\noindent{\bf Problem 5.} Find suitable conditions under which each torsion-free ISI-group of an idempotent type is an IFI-group.

\medskip

In a way of similarity to Theorem~\ref{square} and Corollary~\ref{cartesian}, both listed above, we may state our next query.

\medskip

\noindent{\bf Problem 6.} For a group $G$ of an idempotent type, does it follow that the square $G\oplus G$ is an ISI-group if, and only if, $G\oplus G$ is an IFI-group?

\medskip

The restriction on the idempotent type is essential and cannot be ignored: in fact, if $G$ is a torsion-free separable group of a non-idempotent type, then it is necessarily si-simple and thus $G\oplus G$ is too si-simple and hence an ISI-group but {\it not} an IFI-group, because $G\oplus G$ is still of a non-idempotent type.

\medskip

\noindent{\bf Acknowledgement.} The authors would like to thank the unknown expert referee for his/her valuable comments and suggestions which led to an improvement of the article's shape.

\medskip

\noindent{\bf Funding:} The scientific work of Andrey R. Chekhlov was supported by the Ministry of Science and Higher Education of Russia (agreement No. 075-02-2023-943). The scientific work of Peter V. Danchev was supported in part by the Bulgarian National Science Fund under Grant KP-06 No 32/1 of December 07, 2019, as well as by the Junta de Andaluc\'ia, Grant FQM 264, and by the BIDEB 2221 of T\"UB\'ITAK.


\bigskip

\begin{thebibliography}{99}

\bibitem{GC}
G. C\u{a}lug\u{a}reanu, {\it Strongly invariant subgroups}, Glasg. Math. J. (2) {\bf 57} (2015), 431--443.

\bibitem{C2}
A. R. Chekhlov, {\it On a direct sum of irredicible groups}, Math. Notes (5) {\bf 97} (2015), 815--817.

\bibitem{C1}
A. R. Chekhlov, {\it On strongly invariant subgroups of abelian groups}, Math. Notes (1) {\bf 102} (2017), 106--110.

\bibitem{CD1}
A. R. Chekhlov and P. V. Danchev, {\it On abelian groups having all proper fully invariant subgroups isomorphic}, Commun. Algebra (12) {\bf 43} (2015), 5059--5073.

\bibitem{CD2}
A. R. Chekhlov and P. V. Danchev, {\it The strongly invariant extending property for abelian groups}, Quaest. Math. (8) {\bf 42} (2019), 997--1017.

\bibitem{CD3}
A. R. Chekhlov and P. V. Danchev, {\it On the socles of strongly inert subgroups of abelian $p$-groups}, Siber. Math. J. (2) {\bf 64} (2023), 459--468.

\bibitem{CD4}
A. R. Chekhlov and P. V. Danchev, {\it On abelian groups having isomorphic proper characteristic subgroups}, J. Commut. Algebra (4) {\bf 15} (2023).

\bibitem{Fo}
A. A. Fomin, {\it Abelian groups with free subgroups of infinite index and their endomorphism groups}, Math. Notes (2) {\bf 36} (1984), 581-–585.

\bibitem{F}
L. Fuchs, {\it Infinite Abelian Groups}, Vol. {\bf I} and Vol. {\bf II}, Academic Press, New York and London, 1970 and 1973.

\bibitem{F1}
L. Fuchs, {\it Abelian Groups}, Springer (Cham, Switzerland, 2015).

\bibitem{GN1}
S. Ya. Grinshpon and M. M. Nikolskaya, {\it Primary IF-groups}, Tomsk State Univ. J. Math. \& Mech. (3) {\bf 15} (2011), 25--31.

\bibitem{GN2}
S. Ya. Grinshpon and M. M. Nikolskaya (Savinkova), {\it Fully invariant subgroups of Abelian $p$-groups with finite Ulm-Kaplansky invariants}, Commun. Algebra (11) {\bf 39} (2011), 4273--4282.

\bibitem{GN3}
S. Ya. Grinshpon and M. M. Nikolskaya, {\it Proper fully invariant subgroups of a torsion-free group isomorphic to the group}, Tomsk State Univ. J. Math. \& Mech. (1) {\bf 17} (2012), 11--15.

\bibitem{Kap}
I. Kaplansky, {\it Infinite Abelian Groups}, University of Michigan Press, Ann Arbor, 1954 and 1969.

\bibitem{KMT}
P. A. Krylov, A. V. Mikhalev and A. A. Tuganbaev, Endomorphism Rings of Abelian Groups, Dordrecht, Boston and London: Kluwer Acad. Publ. (2003).

\end{thebibliography}
\end{document}